\documentclass[11pt]{article}
\usepackage[letterpaper,margin=1in]{geometry}
\usepackage[T1]{fontenc}
\usepackage{lmodern}
\usepackage{microtype}
\usepackage{amsmath}
\usepackage{amssymb}
\usepackage{amsthm}
\usepackage{booktabs}
\usepackage{graphicx}
\usepackage{xcolor}
\usepackage{tikz}
\usetikzlibrary{arrows.meta}
\usepackage[caption=false,font=footnotesize,labelformat=simple]{subfig}%
\usepackage{cite}
\usepackage[colorlinks=true,linkcolor={red!60!black},citecolor={green!45!black},urlcolor={red!60!black}]{hyperref}

\theoremstyle{plain}
\newtheorem{thm}{Theorem}
\newtheorem{lem}{Lemma}
\newtheorem{assum}{Assumption}
\theoremstyle{remark}
\newtheorem{rem}{Remark}

\title{Robust Decoupled Input Filtering for Internal Model Principle-Based Output Regulation of Polynomial Exogenous Signals\thanks{This work was supported in part by the National Science Foundation
under Awards No.~2427222 and No.~2536192, and by the Office of Naval Research
under Award No.~N00014-24-1-2426.}}

\author{Venanzio Cichella\thanks{Department of Mechanical Engineering, University of Iowa, Iowa City, IA 52242, USA (\texttt{venanzio-cichella@uiowa.edu}).}
\and
Lorenzo Marconi\thanks{Department of Electrical, Electronic, and Information Engineering (DEI), University of Bologna, 40136 Bologna, Italy (\texttt{lorenzo.marconi@unibo.it}).}}

\date{}

\begin{document}

\maketitle

\begin{abstract}
Internal model principle-based controllers regulate polynomial references and
disturbances exactly, but only asymptotically, and their transient input can
peak and oscillate, increasingly so as the degree of the exogenous signals
grows. A filter placed in series with the plant smooths the command but alters
the closed-loop spectrum and can destroy stability. We propose an
architecture that keeps the filter outside the loop the regulator stabilizes: the internal model is
driven by the error of an ideal, unfiltered copy of the plant, the physical
plant receives the filtered command, and an auxiliary plant copy absorbs the
mismatch between the two inputs. If the filter satisfies blocking conditions
at the origin, any gains that stabilize the unfiltered loop can be combined
with any such filter without losing regulation of the true error. Under
additive plant uncertainty, we derive a small-gain condition that guarantees
robust stability and regulation.  On a scalar plant, the proposed controller
produces a smoother input and a larger input-delay margin than the unfiltered
one, and it maintains regulation under a large loss of actuator effectiveness.
On a three-mass chain, its tolerance to unmodeled actuator dynamics and delay
stays essentially constant as the regulator is made more aggressive.
\end{abstract}


\section{Introduction}

Output regulation asks that the error between a plant output
and a reference be driven to zero in spite of exogenous disturbances. When
references and disturbances are generated by a known exosystem, the internal
model principle
\cite{francis1975internal,francis1977linear,davison1976robust,wonham1974linear}
states that robust regulation requires the controller to embed a copy of the
exosystem dynamics. The principle has since been extended to nonlinear,
adaptive and hybrid settings
\cite{isidori1990output,huang1990nonlinear,byrnes1997structurally,serrani2001semiglobal,byrnes2003limit,huang2004nonlinear,cox2016isolating},
and it remains the organizing idea behind a large body of regulation schemes
\cite{byrnes2000output,bin2022internal}.

This paper is concerned with polynomial exogenous signals. Our
interest in this class stems from its ability to represent signals of
practical relevance. Polynomials in time, and piecewise-polynomial signals
such as B\'ezier and Bernstein curves, are a standard parameterization of
trajectories in trajectory generation and motion planning
\cite{cichella2020optimal,kielas2022bernstein,maclin2024optimal,farouki2012bernstein},
and, by the Weierstrass approximation theorem, polynomials approximate any
continuous signal uniformly on a compact time interval
\cite{powell1981approximation}. Robust linear output regulation of
B\'ezier-based exogenous signals was first addressed in
\cite{marconi2025robust}. The approach was then extended to piecewise
Bernstein/B\'ezier exogenous signals \cite{cichella2026robust} and to
polynomial signals with time-varying coefficients, namely B\'ezier
trajectories with time-varying control points \cite{mimmo2026robust}.
Polynomial exosystems have also been used as function approximators for
unstructured exogenous signals, which led to polynomial internal models for
finite-horizon robust output regulation \cite{cichella2026polynomial}, and as
the basis of model predictive planning for autonomous systems
\cite{cichella2026model}. For a polynomial of degree $N$, the internal model
consists of one copy of an $(N+1)$-dimensional nilpotent signal generator per
error channel, the simplest realization being a chain of $N+1$ integrators.

The internal model, however, inherits its complexity from the exosystem,
whatever the class of exogenous signals considered. The controller must embed
a copy of every exosystem mode, replicated for each error channel, so the
order of the internal model grows with the dimension of the exosystem and with
the number of regulated outputs. This order can become very large when the
signals to be tracked or rejected are rich, for instance multi-tone or
high-order generators.

Moreover, the internal model principle guarantees only convergence of the
regulated error; it says little about how the error, and in particular the
control input, behave along the way. The transient is governed by the
closed-loop matrix and by the initial mismatch with respect to the steady
state, and both deteriorate as the exosystem order grows: each additional mode
adds a state to the internal model, and keeping the regulation speed
independent of that order requires fast pole placement. The resulting input
can exhibit significant peaking and high-frequency oscillations, a phenomenon
well known in high-gain designs \cite{sussmann1991peaking}. In practice, such
inputs may exceed actuator bandwidth and amplitude limits and excite unmodeled
dynamics.

The obvious remedy is to low-pass filter the command before it reaches the
actuator. A filter inserted in series with the plant, however, becomes part of
the feedback loop. It modifies the closed-loop spectrum that the gains were
designed to place, erodes stability margins, and may destabilize the loop
altogether \cite{skogestad2005multivariable}. Designing the regulator and
filter independently can therefore lead to destabilizing interactions,
regardless of the internal model order. Avoiding such interactions requires a
joint design of the filter and regulator, sacrificing the simplicity of the
internal model principle-based approach.

The architecture proposed here avoids this trade-off by letting the filter
shape the plant input without entering the loop that the regulator stabilizes.
Its starting point is $\mathcal{L}_1$ adaptive control
\cite{hovakimyan2010l1}, where a low-pass filter in the control channel
decouples the speed of adaptation from the robustness of the closed loop, so
that fast adaptation need not produce a high-frequency control signal. We
transfer this filtering structure to internal model principle-based
regulation. The internal model is driven by the error of an \emph{ideal},
unfiltered copy of the plant. The physical plant receives the filtered
command, and an auxiliary state, driven by the mismatch between the filtered
and the unfiltered inputs, removes the effect of the filter from the measured
state. Using a plant model to take an input-side element out of the feedback
loop is itself a familiar device: model-recovery anti-windup
\cite{teel1997l2,zaccarian2011modern} does so for a saturation, and internal
model control and the Smith predictor \cite{morari1989robust,smith1957closer}
do so for a delay. Here the element is a linear filter that we insert
deliberately, and output regulation requires one further ingredient, a set of
blocking conditions at the origin. These conditions ensure that the auxiliary
state vanishes along the steady state, so that regulating the ideal error is
enough to regulate the true error.

A preliminary version of this work has been submitted to the
2027 European Control Conference \cite{cichella2027decoupled}. That paper
considers exactly known plants and polynomial exogenous signals, and its
theoretical result is the decoupling property of Lemma~\ref{lem:convergence}.
The present paper extends the analysis to plants with
additive uncertainty, for which a small-gain condition guarantees robust
stability and regulation (Theorem~\ref{thm:robust}). The numerical study is
extended accordingly to the robustness properties of the proposed
architecture.

The paper makes three contributions.
\begin{itemize}
\item A filtered architecture for internal model principle-based control that
preserves internal stability and asymptotic regulation for \emph{any} gains
that stabilize the unfiltered loop, combined with \emph{any} stable filter
satisfying the blocking conditions
(Lemma~\ref{lem:convergence}, first reported in \cite{cichella2027decoupled}). Control design and filter
design are therefore decoupled.
\item A robustness analysis for additive uncertainty in the plant matrices.
The uncertainty closes a feedback loop around the nominal cascade, and a
small-gain condition on this loop preserves robust stability and regulation
(Theorem~\ref{thm:robust}).
\item Numerical evidence on a scalar plant and on a lightly damped three-mass
chain, comparing the proposed architecture with the unfiltered internal model
principle-based controller and with the naive series filter. On the scalar
plant, the proposed architecture yields smoother inputs, input-delay margins
that are several times larger (in some cases by more than an order of
magnitude), and certified robust regulation under a $55\%$ loss of actuator
effectiveness, where the unfiltered controller fails. On the mass chain, its
tolerance to unmodeled actuator dynamics and delay is essentially independent
of how aggressively the regulator is tuned. The larger margins
relative to the unfiltered design are a numerical observation, and their
analytical characterization is left to future work.
\end{itemize}

Section~\ref{sec:problem} formulates the problem and examines the naive
filtering approach. Section~\ref{sec:proposed} introduces the proposed
architecture and establishes its convergence properties.
Section~\ref{sec:robust} extends the analysis to plant
uncertainty,
Section~\ref{sec:numerics} reports the numerical examples, and
Section~\ref{sec:conclusions} concludes the paper.

\section{Problem Formulation}\label{sec:problem}

Consider the linear system
\begin{align}
\dot{x} &= Ax + Bu + Pw, \label{eq:plant}\\
e &= Cx + Qw, \label{eq:err}
\end{align}
where $x(t)\in\mathbb{R}^n$ is the state, $u(t)\in\mathbb{R}^m$ the control
input, $e(t)\in\mathbb{R}^p$ the regulated error, and $w(t)\in\mathbb{R}^{N+1}$
a basis of the polynomials of degree at most $N$. Since the derivative of a
polynomial of degree at most $N$ is again such a polynomial, $\dot w$ is a
linear combination of the entries of $w$,
\begin{equation}\label{eq:exo}
\dot w = D_N w,
\end{equation}
where $D_N\in\mathbb{R}^{(N+1)\times(N+1)}$ is the differentiation matrix of
the chosen basis, nilpotent of
order $N+1$, with characteristic polynomial $\lambda^{N+1}$ and
$\sigma(D_N)=\{0\}$; \eqref{eq:exo} plays the role of the exosystem. The
matrices are $A\in\mathbb{R}^{n\times n}$,
$B\in\mathbb{R}^{n\times m}$, $C\in\mathbb{R}^{p\times n}$,
$P\in\mathbb{R}^{n\times(N+1)}$ and $Q\in\mathbb{R}^{p\times(N+1)}$, and the
columns of $P$ and $Q$ are the coefficients of the polynomial disturbance and
reference. The objective is to design $u$ so that the closed loop is
internally stable and $e(t)\to 0$.

Throughout, we make the following standing assumption.
\begin{assum}\label{asm:plant}
The pair $(A,B)$ is stabilizable, and
\begin{equation}\label{eq:rank}
\operatorname{rank}\begin{bmatrix} A & B\\ C & 0\end{bmatrix}=n+p.
\end{equation}
\end{assum}
\begin{sloppypar}
Condition \eqref{eq:rank} is the non-resonance condition
$\operatorname{rank}\big[\begin{smallmatrix}\lambda I-A & B\\ C & 0\end{smallmatrix}\big]=n+p$
of linear output regulation \cite{davison1976robust,francis1977linear},
evaluated at the only eigenvalue $\lambda=0$ of the exosystem \eqref{eq:exo}.
It states that the plant has no transmission zero at the origin, and it
requires $p\le m$. Together with stabilizability, it is necessary and
sufficient for the existence of gains that stabilize the plant augmented with
the internal model introduced below \cite{davison1976robust}.
\end{sloppypar}

For simplicity, we also take $A$ Hurwitz, which entails no loss of generality.
Since $(A,B)$ is stabilizable, a preliminary state feedback $u = K_0 x + v$,
with $K_0\in\mathbb{R}^{m\times n}$ and $v\in\mathbb{R}^m$, can always place
the modes of $A+BK_0$ at desired stable locations; the design then proceeds on
the prestabilized plant, with $A+BK_0$ in place of $A$ and $v$ as the new
input. 

For polynomial exogenous signals, we adopt the internal model unit of
\cite{cichella2026polynomial}:
\begin{align}
\dot{\eta} &= D\eta + \Upsilon e, \label{eq:im}\\
u &= K_x x + K_\eta \eta,
\label{eq:imcontrol}
\end{align}
with $\eta\in\mathbb{R}^{n_\eta}$, $n_\eta:=(N+1)p$, gains
$K_x\in\mathbb{R}^{m\times n}$, $K_\eta\in\mathbb{R}^{m\times n_\eta}$, and
\begin{equation}\label{eq:DUps}
D=I_p\otimes D_N\in\mathbb{R}^{n_\eta\times n_\eta},\qquad
\Upsilon=I_p\otimes\Upsilon_N\in\mathbb{R}^{n_\eta\times p}.
\end{equation}
Here $D_N$ is the differentiation matrix in \eqref{eq:exo}, and
$\Upsilon_N\in\mathbb{R}^{N+1}$ is any vector such that the pair
$(D_N,\Upsilon_N)$ is controllable, in which case $(D,\Upsilon)$ is
controllable as well \cite{cichella2026polynomial}. The internal model thus
contains one copy of the exosystem \eqref{eq:exo} per error channel, as the
internal model principle requires.

We make explicit why this structure enforces regulation, since the same
argument is used repeatedly below. Differentiating \eqref{eq:im}
$N+1$ times and using $D^{N+1}=0$ gives
\begin{equation}\label{eq:eta_derivative}
\eta^{(N+1)}=\sum_{i=0}^{N}D^{N-i}\Upsilon\,e^{(i)}
=\mathcal{C}(D,\Upsilon)
\begin{bmatrix} e^{(N)}\\ \vdots\\ e\end{bmatrix},
\end{equation}
where $\mathcal{C}(D,\Upsilon)=[\,\Upsilon\ \ D\Upsilon\ \cdots\ D^N\Upsilon\,]\in\mathbb{R}^{n_\eta\times n_\eta}$
is the controllability matrix of $(D,\Upsilon)$. If $\eta$ is a polynomial of
degree at most $N$, the left-hand side vanishes. Since $(D,\Upsilon)$ is
controllable, $\mathcal{C}(D,\Upsilon)$ is square and nonsingular, so $e$ and
its first $N$ derivatives must vanish, that is, $e\equiv0$
\cite{cichella2026polynomial}. The argument uses only $D^{N+1}=0$ and the
controllability of $(D,\Upsilon)$. Consequently, $D_N$ in \eqref{eq:DUps} may
be replaced by the differentiation matrix $\hat D_N$ of any other polynomial
basis of degree $N$, together with a vector $\Upsilon_N$ that makes
$(\hat D_N,\Upsilon_N)$ controllable. Two such bases are related by an
invertible change of coordinates, so $\hat D_N$ is similar to $D_N$ and is
again nilpotent of order $N+1$. 

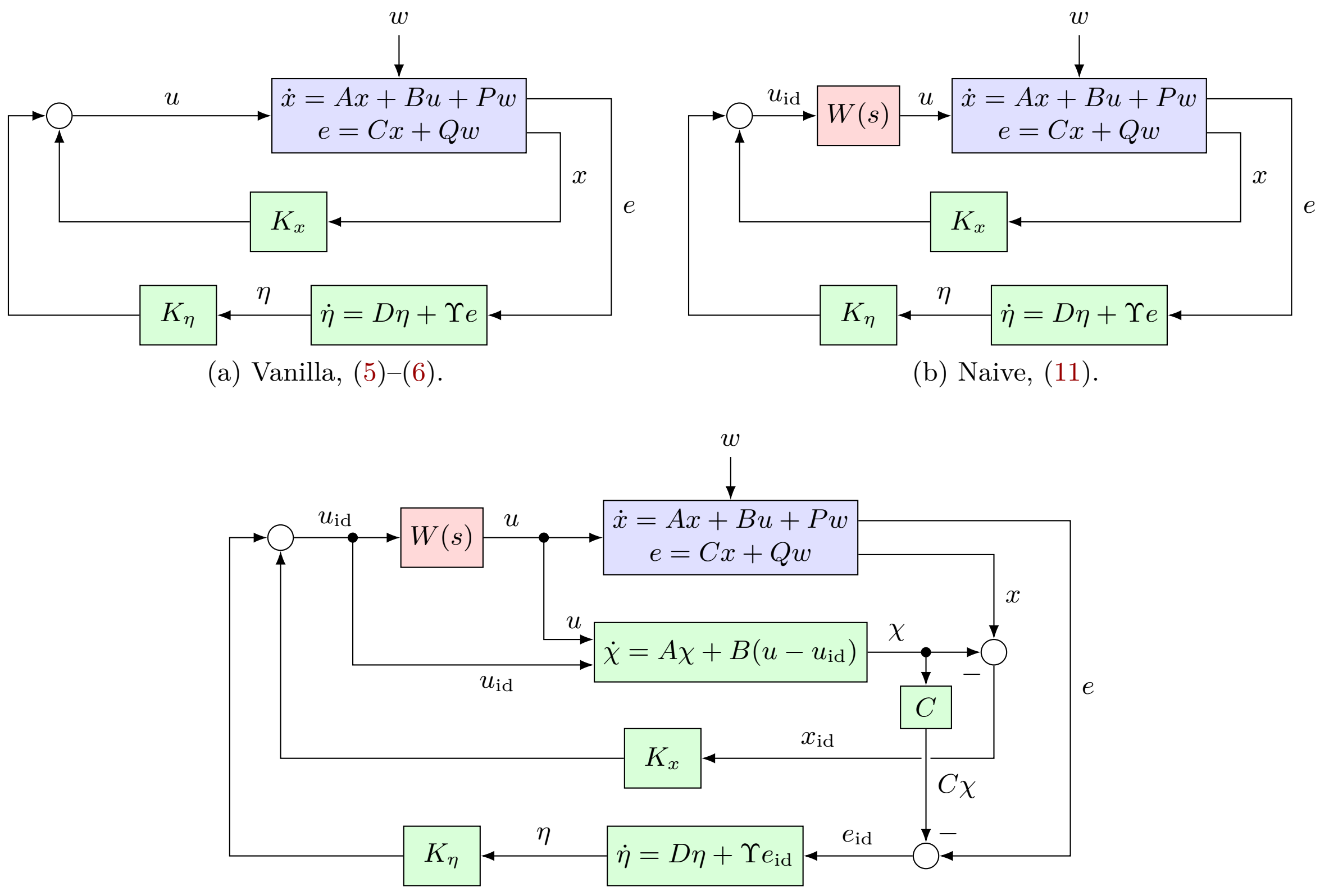
\begin{figure*}[t]
\centering
\tikzset{
  imfblk/.style={draw, rectangle, minimum height=0.7cm, minimum width=0.9cm,
    inner sep=3pt, align=center},
  imfsum/.style={draw, circle, inner sep=0pt, minimum size=0.3cm},
  imfsig/.style={-{Latex[length=1.6mm,width=1.2mm]}},
  imftap/.style={circle, fill, inner sep=0pt, minimum size=1.2mm},
  imfsys/.style={fill=blue!12},
  imfctrl/.style={fill=green!15},
  imffilt/.style={fill=red!15},
}
\subfloat[Vanilla, \eqref{eq:im}--\eqref{eq:imcontrol}.]{%
\begin{tikzpicture}[font=\footnotesize]
\node[imfsum] (s1) at (0,0) {};
\node[imfblk, imfsys, minimum width=2.6cm] (pl) at (4.0,0) {$\dot x=Ax+Bu+Pw$\\ $e=Cx+Qw$};
\node[imfblk, imfctrl] (kx) at (2.7,-1.25) {$K_x$};
\node[imfblk, imfctrl] (im) at (4.0,-2.35) {$\dot\eta=D\eta+\Upsilon e$};
\node[imfblk, imfctrl] (ke) at (1.4,-2.35) {$K_\eta$};
\draw[imfsig] (s1) -- node[above]{$u$} (pl);
\draw[imfsig] (4.0,0.95) node[above]{$w$} -- (pl.north);
\draw[imfsig] ([yshift=-2mm]pl.east) -- (5.9,-0.2) -- node[right]{$x$} (5.9,-1.25) -- (kx);
\draw[imfsig] ([yshift=2mm]pl.east) -- (6.5,0.2) -- node[right]{$e$} (6.5,-2.35) -- (im);
\draw[imfsig] (im) -- node[above]{$\eta$} (ke);
\draw[imfsig] (kx) -- (0,-1.25) -- (s1);
\draw[imfsig] (ke) -- (-0.6,-2.35) -- (-0.6,0) -- (s1);
\end{tikzpicture}%
\label{fig:arch_vanilla}}
\hfil
\subfloat[Naive, \eqref{eq:naive}.]{%
\begin{tikzpicture}[font=\footnotesize]
\node[imfsum] (s1) at (0,0) {};
\node[imfblk, imffilt] (w) at (1.4,0) {$W(s)$};
\node[imfblk, imfsys, minimum width=2.6cm] (pl) at (4.0,0) {$\dot x=Ax+Bu+Pw$\\ $e=Cx+Qw$};
\node[imfblk, imfctrl] (kx) at (2.7,-1.25) {$K_x$};
\node[imfblk, imfctrl] (im) at (4.0,-2.35) {$\dot\eta=D\eta+\Upsilon e$};
\node[imfblk, imfctrl] (ke) at (1.4,-2.35) {$K_\eta$};
\draw[imfsig] (s1) -- node[above]{$u_{\text{id}}$} (w);
\draw[imfsig] (w) -- node[above]{$u$} (pl);
\draw[imfsig] (4.0,0.95) node[above]{$w$} -- (pl.north);
\draw[imfsig] ([yshift=-2mm]pl.east) -- (5.9,-0.2) -- node[right]{$x$} (5.9,-1.25) -- (kx);
\draw[imfsig] ([yshift=2mm]pl.east) -- (6.5,0.2) -- node[right]{$e$} (6.5,-2.35) -- (im);
\draw[imfsig] (im) -- node[above]{$\eta$} (ke);
\draw[imfsig] (kx) -- (0,-1.25) -- (s1);
\draw[imfsig] (ke) -- (-0.6,-2.35) -- (-0.6,0) -- (s1);
\end{tikzpicture}%
\label{fig:arch_naive}}
\\
\subfloat[Proposed, \eqref{eq:internal_model}--\eqref{eq:control_filter}.]{%
\begin{tikzpicture}[font=\footnotesize]
\node[imfsum] (s1) at (0,0) {};
\node[imfblk, imffilt] (w) at (1.9,0) {$W(s)$};
\node[imfblk, imfsys, minimum width=2.6cm] (pl) at (5.3,0) {$\dot x=Ax+Bu+Pw$\\ $e=Cx+Qw$};
\node[imfblk, imfctrl] (cp) at (5.3,-1.35) {$\dot\chi=A\chi+B(u-u_{\text{id}})$};
\node[imfsum] (s3) at (8.4,-1.35) {};
\node[imfblk, imfctrl, minimum width=0.6cm, minimum height=0.5cm] (c) at (7.6,-2.0) {$C$};
\node[imfsum] (s4) at (7.6,-3.75) {};
\node[imfblk, imfctrl] (kx) at (4.5,-2.6) {$K_x$};
\node[imfblk, imfctrl] (im) at (5.0,-3.75) {$\dot\eta=D\eta+\Upsilon e_{\text{id}}$};
\node[imfblk, imfctrl] (ke) at (1.9,-3.75) {$K_\eta$};
\coordinate (ta) at (0.85,0);
\coordinate (tb) at (3.1,0);
\coordinate (tc) at (7.6,-1.35);
\draw[imfsig] (s1) -- (w);
\node[above] at (0.65,0) {$u_{\text{id}}$};
\draw[imfsig] (w) -- (pl);
\node[above] at (2.72,0) {$u$};
\draw[imfsig] (5.3,0.95) node[above]{$w$} -- (pl.north);
\draw[imfsig] (tb) |- node[pos=0.8, above]{$u$} ([yshift=1.5mm]cp.west);
\draw[imfsig] (ta) |- node[pos=0.8, below]{$u_{\text{id}}$} ([yshift=-1.5mm]cp.west);
\draw (cp) -- node[above]{$\chi$} (tc);
\draw[imfsig] (tc) -- (s3);
\node at ([xshift=-2.6mm, yshift=-2.6mm]s3) {$-$};
\draw[imfsig] ([yshift=-2mm]pl.east) -- (8.4,-0.2) -- node[right]{$x$} (s3);
\draw[imfsig] (s3) -- (8.4,-2.6) -- node[above, pos=0.6]{$x_{\text{id}}$} (kx);
\draw[imfsig] (tc) -- (c);
\draw[imfsig, preaction={draw=white, line width=3pt, -}] (c) -- node[pos=0.5, right]{$C\chi$} (s4);
\node at ([xshift=2.6mm, yshift=2.6mm]s4) {$-$};
\draw[imfsig] ([yshift=2mm]pl.east) -- (9.3,0.2) -- node[right]{$e$} (9.3,-3.75) -- (s4);
\draw[imfsig] (s4) -- node[above]{$e_{\text{id}}$} (im);
\draw[imfsig] (im) -- node[above]{$\eta$} (ke);
\draw[imfsig] (kx) -- (0,-2.6) -- (s1);
\draw[imfsig] (ke) -- (-0.6,-3.75) -- (-0.6,0) -- (s1);
\node[imftap] at (ta) {};
\node[imftap] at (tb) {};
\node[imftap] at (tc) {};
\end{tikzpicture}%
\label{fig:arch_proposed}}
\caption{Block diagrams of the three architectures, with the
plant \eqref{eq:plant}--\eqref{eq:err} driven by the exogenous signal $w$.
(a)~Vanilla internal model principle-based controller. (b)~Naive series filter.
(c)~Proposed architecture, with ideal-plant state $x_{\text{id}}=x-\chi$ and
ideal error $e_{\text{id}}=e-C\chi$, as in
\eqref{eq:ideal_plant}--\eqref{eq:ideal_error}. Under the uncertainty
\eqref{eq:uncertain}, the plant copy uses $(A_{\mathrm n},B_{\mathrm n})$, as in
\eqref{eq:obs_r}.}
\label{fig:arch}
\end{figure*}

Substituting \eqref{eq:im}--\eqref{eq:imcontrol} into
\eqref{eq:plant}--\eqref{eq:err} gives the closed loop, in the state
$z=[\,x^\top\ \eta^\top\,]^\top\in\mathbb{R}^{n+n_\eta}$,
\begin{equation} \label{eq:cl}
\begin{bmatrix}\dot{x}\\ \dot{\eta}\end{bmatrix}
=\underbrace{\begin{bmatrix} A+BK_x & BK_\eta \\ \Upsilon C & D \end{bmatrix}}_{F}
\begin{bmatrix}x\\ \eta\end{bmatrix}
+\underbrace{\begin{bmatrix} P \\ \Upsilon Q\end{bmatrix}}_{G} w ,
\end{equation}
with $F\in\mathbb{R}^{(n+n_\eta)\times(n+n_\eta)}$ and
$G\in\mathbb{R}^{(n+n_\eta)\times(N+1)}$.
Fig.~\ref{fig:arch_vanilla} depicts this loop.

Under Assumption~\ref{asm:plant}, and since $(D,\Upsilon)$ is controllable,
there exist gains $K_x,K_\eta$ such that $F$ is Hurwitz. Let
$z^\star(t)=[\,x^\star(t)^\top\ \eta^\star(t)^\top\,]^\top$ denote the polynomial steady-state
solution of \eqref{eq:cl} corresponding to $w(t)$. By the internal model
principle \cite{francis1975internal}, the internal model in \eqref{eq:im}
guarantees that the error vanishes identically on this solution,
\[
e^\star(t)=Cx^\star(t)+Qw(t)\equiv0.
\]
Since $F$ is Hurwitz, the deviation $\tilde z(t)=z(t)-z^\star(t)$ satisfies the
homogeneous system
\[
\dot{\tilde z}=F\tilde z,
\]
and therefore
\[
\tilde z(t)=e^{Ft}\tilde z(0).
\]
Consequently,
\begin{equation} \label{eq:etransient}
e(t)=C\Gamma\tilde z(t)
=C\Gamma e^{Ft}\bigl(z(0)-z^\star(0)\bigr),
\end{equation}
where $\Gamma=[\,I_n\;\;0\,]\in\mathbb{R}^{n\times(n+n_\eta)}$ selects the
plant state from $z$. Hence $e(t)\to0$
exponentially, at a rate determined solely by $F$ and independent of the
coefficients of the polynomial disturbance and reference.

Equation~\eqref{eq:etransient} also exposes the limits of this guarantee.
Regulation is exact only asymptotically, and the transient error depends on
both the matrix $F$ and the initial mismatch $z(0)-z^\star(0)$. Increasing $N$
generally worsens the transient in both respects. Most importantly for this
paper, the corresponding control input can exhibit significant peaking and
oscillatory behavior during the transient, and these effects become more
pronounced as $N$ increases.

The most direct way to mitigate this behavior is to filter the control input.
The internal model principle-based controller \eqref{eq:im} is left unchanged, but its output
is treated as an \emph{ideal} command $u_{\text{id}}(t)\in\mathbb{R}^m$ and
passed through a filter before reaching the plant
(Fig.~\ref{fig:arch_naive}):
\begin{equation}
\label{eq:naive}
\begin{aligned}
\dot{\eta} &= D\eta + \Upsilon e, \\
u_{\text{id}} &= K_x x + K_\eta \eta, \\
u(s) &= W(s)\,u_{\text{id}}(s), \end{aligned}
\end{equation}
where $W(s)$ is an $m\times m$ strictly proper filter chosen to attenuate the high-frequency
content and peaking of the raw input $u_{\text{id}}$. As is well known from
classical control theory, however, inserting such a filter in the loop can
compromise closed-loop stability, erode robustness margins, and degrade
overall performance.

Our aim is a filtering architecture that attenuates these undesirable input
transients while preserving the closed-loop performance of the system.

\section{Proposed approach}\label{sec:proposed}

We consider the controller, depicted in
Fig.~\ref{fig:arch_proposed},
\begin{align}
    \dot{\eta}(t) &= D\eta(t) + \Upsilon\big(e(t) - C\chi(t)\big), \label{eq:internal_model}\\
    \dot{\chi}(t) &= A\chi(t) + B\big(u(t) - u_{\text{id}}(t)\big), \label{eq:observer}\\
    u_{\text{id}}(t) &= K_x\big(x(t) - \chi(t)\big) + K_{\eta}\eta(t), \label{eq:ideal_control}\\
    u(s) &= W(s)\,u_{\text{id}}(s). \label{eq:control_filter}
\end{align}
Here $u_{\text{id}}$ is the state-feedback law that would be applied in the
absence of filtering, $\chi(t)\in\mathbb{R}^n$ is an internal copy of the plant
driven by the input mismatch $u-u_{\text{id}}$, and $u$ is obtained by passing $u_{\text{id}}$
through the $m\times m$ strictly proper filter $W(s)$, which satisfies
\begin{equation}
W(0)=I_m,\qquad \frac{d^k W}{ds^k}(0)=0,\quad k=1,\dots,N. \label{eq:blocking}
\end{equation}
The filter is stable, and we fix a minimal realization
\begin{equation}\label{eq:Wreal}
\begin{aligned}
& \dot\xi = A_W\xi + B_W u_{\text{id}}, \\
& u = C_W\xi,\\
& W(s) = C_W(sI-A_W)^{-1}B_W,
\end{aligned}
\end{equation}
with $\xi\in\mathbb{R}^{n_W}$, $A_W\in\mathbb{R}^{n_W\times n_W}$ Hurwitz,
$B_W\in\mathbb{R}^{n_W\times m}$, $C_W\in\mathbb{R}^{m\times n_W}$, and
$W(0)=-C_W A_W^{-1}B_W = I_m$.

The reasoning behind this structure is clearest by comparison with the naive
scheme of Equation~\eqref{eq:naive}
(compare Figs.~\ref{fig:arch_naive}
and~\ref{fig:arch_proposed}). When the input is filtered directly, the
filter lag propagates into the error $e$ that drives the internal model
\eqref{eq:naive}: $e$ is the true error of the plant \eqref{eq:plant}, whose
input has been delayed by $W(s)$. The internal model is thus forced to run at
the (slow) speed of the filtered loop, and this is the source of the stability
and robustness degradation noted above.

The proposed architecture separates these two time scales. The physical plant
is driven by the filtered input $u(s)=W(s)u_{\text{id}}(s)$, so the actuator
sees a smooth command. The internal model, in contrast, is driven by the error
of an \emph{ideal}, unfiltered copy of the plant, and is therefore free to adapt
quickly. To make this precise, define the ideal-plant state
$x_{\text{id}} = x - \chi\in\mathbb{R}^n$. Differentiating and using \eqref{eq:plant} and
\eqref{eq:observer} gives
\begin{equation}
\dot{{x}}_{\text{id}}
= A{x}_{\text{id}} + Bu_{\text{id}} + Pw,
\label{eq:ideal_plant}
\end{equation}
so ${x}_{\text{id}}$ is the state the plant \emph{would} have if it were driven
by the unfiltered input $u_{\text{id}}$; the filter has been cancelled out. The
signal feeding the internal model in \eqref{eq:internal_model} is precisely the
error of this ideal plant,
\begin{equation}
e - C\chi = C(x-\chi) + Qw = C{x}_{\text{id}} + Qw =: {e}_{\text{id}},
\label{eq:ideal_error}
\end{equation}
where $e_{\text{id}}\in\mathbb{R}^p$ is unaffected by $W(s)$. The internal
model therefore regulates
${e}_{\text{id}}$ on the fast, unfiltered dynamics \eqref{eq:ideal_plant}, while
the command $u$ delivered to the physical system remains smooth. In effect,
the two loops run at different speeds: fast adaptation on the ideal plant,
gentle actuation on the real one. This separation mirrors the role of the
control filter in $\mathcal{L}_1$ adaptive control \cite{hovakimyan2010l1},
with the internal model in place of the adaptive law.

What remains to be shown is that the closed loop is internally stable and that
regulating $e_{\text{id}}$ also regulates the true error $e$. This is the
content of the following result, on which the rest of the paper builds.
\begin{lem}\label{lem:convergence}
Consider the plant \eqref{eq:plant}--\eqref{eq:err} driven by the controller
\eqref{eq:internal_model}--\eqref{eq:control_filter}, with $w(t)$ a polynomial basis of
degree at most $N$. Suppose that
\begin{enumerate}
    \item[(i)] the gains $K_x,K_\eta$ render
    \begin{equation} \label{eq:clmatrixF}
    F = \begin{bmatrix} A+BK_x & BK_\eta \\ \Upsilon C & D \end{bmatrix}
    \end{equation}
    Hurwitz; and
    \item[(ii)] the filter $W(s)$ satisfies the blocking conditions
    \eqref{eq:blocking}.
\end{enumerate}
Then the closed loop is internally stable and the regulated error
\eqref{eq:err} satisfies $e(t)\to 0$ as $t\to\infty$.
\end{lem}

\begin{proof}
Introduce the ideal-plant state $x_{\text{id}} = x - \chi$ and
stack $z_{\text{id}} = [\,x_{\text{id}}^\top\ \ \eta^\top\,]^\top\in\mathbb{R}^{n+n_\eta}$. By
\eqref{eq:ideal_plant}, and because the internal model \eqref{eq:internal_model}
is driven by $e - C\chi = Cx_{\text{id}} + Qw$, the pair $(x_{\text{id}},\eta)$
evolves autonomously as
\begin{equation}\label{eq:zid_dynamics}
\dot z_{\text{id}} = F z_{\text{id}} + G w,
\qquad
G = \begin{bmatrix} P \\ \Upsilon Q \end{bmatrix},
\end{equation}
with $F$ as in \eqref{eq:clmatrixF}. Since $F$ is Hurwitz and
$\sigma(D_N)=\{0\}$, there is a unique polynomial steady state
$z_{\text{id}}^\star=\Pi w$, with $\Pi\in\mathbb{R}^{(n+n_\eta)\times(N+1)}$
solving $\Pi D_N=F\Pi+G$, of degree at most $N$, and
$z_{\text{id}}-z_{\text{id}}^\star\to0$ exponentially. The internal model
\eqref{eq:internal_model} is \eqref{eq:im} driven by $e_{\text{id}}$, so
\eqref{eq:eta_derivative} holds with $e_{\text{id}}$ in place of $e$. Along the
steady state, $\eta^\star$ is a polynomial of degree at most $N$, so
$(\eta^\star)^{(N+1)}\equiv0$, and the nonsingularity of
$\mathcal{C}(D,\Upsilon)$ gives
\begin{equation}\label{eq:eid_zero}
e_{\text{id}}^\star = C x_{\text{id}}^\star + Q w \equiv 0,
\end{equation}
and $e_{\text{id}}=C(x_{\text{id}}-x_{\text{id}}^\star)\to0$ exponentially.

It remains to account for the filter-induced deviation $\chi$, because the true
error decomposes as
\begin{equation}\label{eq:e_decomp}
e = Cx + Qw = C x_{\text{id}} + Qw + C\chi = e_{\text{id}} + C\chi .
\end{equation}
By \eqref{eq:observer}, $\chi$ is driven by the input mismatch
\begin{equation}\label{eq:chi_dynamics}
\dot\chi = A\chi + B\big(u - u_{\text{id}}\big),
\qquad
u - u_{\text{id}} = \big(W(s) - I_m\big)\,u_{\text{id}} .
\end{equation}
By (i) the ideal loop is internally stable, so its steady-state input
$u_{\text{id}}^\star = K_x x_{\text{id}}^\star + K_\eta \eta^\star$ is a
polynomial in $t$ of degree at most $N$. We treat the steady-state and the
transient parts of the forcing in \eqref{eq:chi_dynamics} separately. The
blocking conditions \eqref{eq:blocking} give $W(s) - I_m$ a zero of order $N+1$
at the origin, so the response of the stable filter $W(s)-I_m$ to the
steady-state input converges to zero,
$\big(W(s)-I_m\big)u_{\text{id}}^\star \to 0$. The transient part
$u_{\text{id}}-u_{\text{id}}^\star
=K_x(x_{\text{id}}-x_{\text{id}}^\star)+K_\eta(\eta-\eta^\star)$ decays
exponentially, and so does its response through $W(s)-I_m$. Hence the forcing
$u-u_{\text{id}}$ in \eqref{eq:chi_dynamics} converges to zero. Since $A$ is
Hurwitz (standing assumption), $\chi$ is the state of an exponentially stable
linear system driven by a vanishing input, and therefore $\chi(t) \to 0$ and
$C\chi(t) \to 0$. Together with \eqref{eq:e_decomp} and $e_{\text{id}}\to 0$,
this yields $e(t) \to 0$.
\end{proof}

\begin{rem}\label{rem:decoupling}
The two hypotheses of Lemma~\ref{lem:convergence} constrain disjoint objects,
so control design and filter design are \emph{decoupled}. Any $K_x,K_\eta$
rendering $F$ in \eqref{eq:clmatrixF} Hurwitz may be combined with any filter
$W(s)$ satisfying the blocking conditions \eqref{eq:blocking}, and regulation is
preserved regardless of the pairing. This is the essential difference from the
naive series filter, where the filter enters the loop that $K_x,K_\eta$ were
designed to stabilize (Fig.~\ref{fig:arch_naive}) and can therefore destroy closed-loop stability. In the
proposed architecture, $W(s)$ acts only on the ideal input and leaves the
spectrum of $F$ untouched, so the filter can be tuned purely for input shaping,
with no feedback interaction with the regulator.
\end{rem}

\section{Partially unknown systems}\label{sec:robust}

Lemma~\ref{lem:convergence} relies on an exact copy of the plant: the
cancellation that produces \eqref{eq:ideal_plant} uses the same $A$ and $B$ in
\eqref{eq:plant} and \eqref{eq:observer}. We now ask what survives when the
plant matrices are known only up to additive uncertainty,
\begin{equation}\label{eq:uncertain}
A = A_{\mathrm{n}} + \Delta A,\qquad B = B_{\mathrm{n}} + \Delta B,
\end{equation}
with $A_{\mathrm{n}}\in\mathbb{R}^{n\times n}$ and
$B_{\mathrm{n}}\in\mathbb{R}^{n\times m}$ used for design, and
$\Delta A\in\mathbb{R}^{n\times n}$ and $\Delta B\in\mathbb{R}^{n\times m}$
unknown. The controller can then only be implemented on the nominal data:
\begin{align}
\dot{\eta} &= D\eta + \Upsilon\big(e - C\chi\big), \label{eq:im_r}\\
\dot{\chi} &= A_{\mathrm{n}}\chi + B_{\mathrm{n}}\big(u - u_{\text{id}}\big), \label{eq:obs_r}\\
u_{\text{id}} &= K_x\big(x - \chi\big) + K_{\eta}\eta, \qquad
u(s) = W(s)\,u_{\text{id}}(s), \label{eq:law_r}
\end{align}
where the filter is realized as in \eqref{eq:Wreal}.

Our objective is a result analogous to Lemma~\ref{lem:convergence} for this
partially unknown system. To state it, we first define
\begin{equation}\label{eq:Fn}
F_{\mathrm n}=\begin{bmatrix} A_{\mathrm n}+B_{\mathrm n}K_x & B_{\mathrm n}K_\eta\\
\Upsilon C & D\end{bmatrix},
\end{equation}
which is the nominal counterpart of the closed-loop matrix $F$ in
\eqref{eq:clmatrixF}, that is, the matrix that would govern the ideal loop in
the absence of uncertainty. We further define the $n\times n$ and $m\times n$
transfer matrices
\begin{align}
T_A(s) &= T_x(s) + (sI-A_{\mathrm n})^{-1}B_{\mathrm n}\big(W(s)-I_m\big)T_u(s), \label{eq:TA}\\
T_B(s) &= W(s)\,T_u(s), \label{eq:TB}
\end{align}
where
\begin{equation}\label{eq:TxTu}
T_x(s)=\Gamma\,(sI-F_{\mathrm n})^{-1}E_z,\qquad
T_u(s)=K\,(sI-F_{\mathrm n})^{-1}E_z,
\end{equation}
with $E_z=[\,I_n\ \ 0\,]^\top\in\mathbb{R}^{(n+n_\eta)\times n}$, $\Gamma$ as in
\eqref{eq:etransient}, and $K=[K_x\ \ K_\eta]$. The role of these maps becomes
clear in the proof; here they are introduced only to state the main
result.

\begin{thm}\label{thm:robust}
Consider the uncertain plant \eqref{eq:uncertain} in closed loop with
\eqref{eq:im_r}--\eqref{eq:law_r}. Suppose (i) $K_x,K_\eta$ render
$F_{\mathrm n}$ in \eqref{eq:Fn} Hurwitz; (ii) $W(s)$ satisfies the
blocking conditions \eqref{eq:blocking}; (iii) $A_{\mathrm n}$ and $A_W$ are
Hurwitz; and (iv) the model uncertainties $\Delta A,\Delta B$ and the filter satisfy the
small-gain condition
\begin{equation}\label{eq:sg}
\big\|\Delta A\,T_A + \Delta B\,T_B\big\|_\infty < 1,
\end{equation}
with $T_A,T_B$ as in \eqref{eq:TA}--\eqref{eq:TB}. Then the closed
loop is internally stable and the regulated error \eqref{eq:err} satisfies
$e(t)\to0$ as $t\to\infty$.
\end{thm}

\begin{proof}
Let $x_{\text{id}}=x-\chi$,
$z_{\text{id}}=[\,x_{\text{id}}^\top\ \eta^\top\,]^\top$, and
$\zeta=[\,z_{\text{id}}^\top\ \chi^\top\ \xi^\top\,]^\top\in\mathbb{R}^{n_\zeta}$,
with $n_\zeta=2n+n_\eta+n_W$. Using \eqref{eq:plant} with
\eqref{eq:uncertain}, together with \eqref{eq:im_r}--\eqref{eq:law_r}
and \eqref{eq:Wreal}, the closed loop reads
\begin{equation}\label{eq:aug}
\dot\zeta = M\zeta
+ \begin{bmatrix}G\\0\\0\end{bmatrix}w,
\end{equation}
with \(M=M_{\mathrm n} + E\,\Delta\,H\),
\begin{equation}\label{eq:EDH}
\begin{aligned}
M_{\mathrm{n}} &=\begin{bmatrix} F_{\mathrm n} & 0 & 0\\ -B_{\mathrm n}K & A_{\mathrm n} & B_{\mathrm n}C_W\\ B_WK & 0 & A_W \end{bmatrix},\qquad
E=\begin{bmatrix} E_z\\0\\0\end{bmatrix},\\
\Delta&=[\,\Delta A\ \ \Delta B\,],\qquad
H=\begin{bmatrix} \Gamma & I_n & 0\\ 0 & 0 & C_W\end{bmatrix},
\end{aligned}
\end{equation}
where $H\zeta=[x^\top\ u^\top]^\top$, since $x=x_{\text{id}}+\chi$ and
$u=C_W\xi$.

We first show that under the assumptions $M:=M_{\mathrm n}+E\Delta H$ is Hurwitz. To this end, we first note that the matrix
$M_{\mathrm n}$ is block lower-triangular with diagonal blocks $F_{\mathrm n}$
and $\big[\begin{smallmatrix}A_{\mathrm n} & B_{\mathrm n}C_W\\ 0 & A_W\end{smallmatrix}\big]$,
so it is Hurwitz by (i) and (iii). Then, we note that for all \(s \notin \sigma(M_n)\) we have 
$$
\begin{aligned}
\det(sI-M)&=\det\big(sI-M_{\mathrm n}-E\Delta H\big)\\
&=\det(sI-M_{\mathrm n})\\
&\qquad\times\det\!\big(I-(sI-M_{\mathrm n})^{-1}E\Delta H\big)\\
&=\det(sI-M_{\mathrm n})\\
&\qquad\times\det\!\big(I_n-\Delta H(sI-M_{\mathrm n})^{-1}E\big)\\
&=\det(sI-M_{\mathrm n})\,\det\!\big(I_n-\Delta T(s)\big).
\end{aligned}
$$
with \(T(s)= H(sI-M_{\mathrm n})^{-1}E\). 
The first factor is nonzero for
$\operatorname{Re}s\ge0$.Since \(M_n\) is Hurwitz, $T$ is analytic and
bounded on $\operatorname{Re}s\ge0$, so by the maximum modulus principle
$\sup_{\operatorname{Re}s\ge0}\|\Delta T(s)\|=\|\Delta T\|_\infty<1$. Hence
$I-\Delta T(s)$ is nonsingular there, and therefore $\det(sI-M)\neq0$ for all $\operatorname{Re}s\ge0$, so $M$ has no eigenvalue in the closed right half-plane, hence is Hurwitz. 

Now split $T$ into its $x$ and $u$ components, $T=[\,T_A^\top\ \ T_B^\top\,]^\top$.
Because $M_{\mathrm n}$ is block lower-triangular, these components can be
computed in cascade. The signal $d$ enters only the $z_{\text{id}}$ block,
whose responses $x_{\text{id}}=T_x\,d$ and $u_{\text{id}}=T_u\,d$ are given by
\eqref{eq:TxTu}. The filter gives $u=W u_{\text{id}}$, and the copy gives
$\chi=(sI-A_{\mathrm n})^{-1}B_{\mathrm n}(W-I_m)u_{\text{id}}$, so that
$x=x_{\text{id}}+\chi$. Hence $T_A$ and $T_B$ are exactly the matrices in
\eqref{eq:TA}--\eqref{eq:TB}. Since $\Delta=[\,\Delta A\ \ \Delta B\,]$, we have
$\Delta T=\Delta A\,T_A+\Delta B\,T_B$, so condition (iv) is equivalent to
$\|\Delta T\|_\infty<1$, which implies that $M$ in \eqref{eq:aug} is Hurwitz.

Finally, we repeat the argument of the proof of Lemma~\ref{lem:convergence}
for \eqref{eq:aug}. Since $M$ is Hurwitz and $\sigma(D_N)=\{0\}$, there is
a unique $\Pi\in\mathbb{R}^{n_\zeta\times(N+1)}$ solving
$\Pi D_N=M\Pi+\big[\begin{smallmatrix}G\\0\\0\end{smallmatrix}\big]$, and
$\zeta-\Pi w\to0$ exponentially. Thus $\zeta$ converges to the polynomial
steady state $\zeta^\star=\Pi w$, of degree at most $N$, and in particular
$\eta$ converges to a polynomial $\eta^\star$ of degree at most $N$. The
internal model \eqref{eq:im_r} is driven by
$e-C\chi=Cx_{\text{id}}+Qw=e_{\text{id}}$, as in \eqref{eq:ideal_error}, so
\eqref{eq:eta_derivative} and the nonsingularity of $\mathcal{C}(D,\Upsilon)$
give $e_{\text{id}}^\star\equiv0$ and hence $e_{\text{id}}\to0$. Likewise,
$u_{\text{id}}=Kz_{\text{id}}$ converges to the polynomial
$u_{\text{id}}^\star=Kz_{\text{id}}^\star$ of degree at most $N$, so the
blocking conditions \eqref{eq:blocking} give $(W-I_m)u_{\text{id}}\to0$, as in
the proof of Lemma~\ref{lem:convergence}. Since $A_{\mathrm n}$ is Hurwitz,
the copy \eqref{eq:obs_r} then gives $\chi\to0$. Therefore, by
\eqref{eq:e_decomp}, $e=e_{\text{id}}+C\chi\to0$.
\end{proof}

\begin{rem}[Role of $W$ in \eqref{eq:sg}]\label{rem:roleW}
Condition \eqref{eq:sg} is sufficient for robust
stability. For $W=I_m$ it reduces to $\|\Delta A\,T_x+\Delta B\,T_u\|_\infty<1$,
which only measures how robust the unfiltered feedback loop is to $\Delta A$ and
$\Delta B$. In the absence of uncertainty, Lemma~\ref{lem:convergence} places no
limit on the filter bandwidth. Under uncertainty,
$\|\Delta A\,T_A+\Delta B\,T_B\|_\infty\le\|\Delta A\|\,\|T_A\|_\infty+\|\Delta B\|\,\|T_B\|_\infty$,
so \eqref{eq:sg} holds whenever the right-hand side is below one. For
$W(s)=W_1(s/\omega_c)$, with $W_1$ a fixed stable, strictly proper filter
satisfying \eqref{eq:blocking}, $\|T_A\|_\infty$ and $\|T_B\|_\infty$ are bounded
uniformly in $\omega_c>0$, because $\|W\|_\infty=\|W_1\|_\infty$. Hence, if
$\|\Delta A\|$ and $\|\Delta B\|$ are small enough, \eqref{eq:sg} holds for every
cutoff $\omega_c$, and the filter can be made as aggressive as desired. For larger uncertainty,
$\omega_c$ must be chosen so that \eqref{eq:sg} holds. Choosing the cutoff that
maximizes smoothing for a given uncertainty bound is an open problem.
\end{rem}

\section{Numerical examples}\label{sec:numerics}

\subsection{Scalar plant}\label{sec:num_setup}

We consider the system in \eqref{eq:plant}--\eqref{eq:err} with $A_{\mathrm n}=-1$, $B_{\mathrm n}=1$, $C=1$ and $Q=0$. The exogenous
signal is the normalized polynomial basis
$w(t)=[\,1,\ t/T,\ \dots,\ (t/T)^N\,]^\top\in\mathbb{R}^{N+1}$ with $T=100$\,s. For each $N$, the
entries of $P\in\mathbb{R}^{1\times(N+1)}$ are
drawn at random (normal, standard deviation $0.2$, fixed seed), so the task is
to reject a random polynomial disturbance of degree $N$ and drive $e = x \to0$.

For the internal model \eqref{eq:im}--\eqref{eq:DUps} we use $D_N$ with ones on the first superdiagonal
together with $\Upsilon_N=[\,0\ \cdots\ 0\ 1\,]^\top$. The gains $K_x,K_\eta$ place the $N+2$ eigenvalues of $F_{\mathrm n}$ in
\eqref{eq:Fn} in a Butterworth pattern on the circle of radius $\rho$, which
fixes the regulation decay rate independently of $N$. For the filter we use
the maximally flat blocking filter
\begin{equation}\label{eq:num_filter}
W(s) = 1 - \Big(\frac{s}{s+\omega_c}\Big)^{N+1},
\end{equation}
realized as a cascade of $N+1$ first-order sections. It is strictly proper and
has all its poles at $-\omega_c$, so $A_W$ is Hurwitz. It also satisfies the
blocking conditions \eqref{eq:blocking} for every $N$ and $\omega_c>0$.

We compare the three architectures discussed in the paper
(Fig.~\ref{fig:arch}):
\emph{vanilla}, the internal model principle-based controller
\eqref{eq:im}--\eqref{eq:imcontrol};
\emph{naive}, the series filter \eqref{eq:naive};
and \emph{proposed}, the controller \eqref{eq:im_r}--\eqref{eq:law_r}, whose
plant copy runs on the nominal data $(A_{\mathrm n},B_{\mathrm n})$ and reduces
to \eqref{eq:internal_model}--\eqref{eq:control_filter} when the model is exact.

The same controllers are tested in two scenarios:
(i) \emph{nominal}: $\Delta A=\Delta B=0$, the setting of
Lemma~\ref{lem:convergence};
(ii) \emph{perturbed}: $\Delta A=0.3$, $\Delta B=-0.55$, i.e.\ the true plant is
$\dot x=-0.7x+0.45u+Pw$. This is the setting of
Theorem~\ref{thm:robust}.

The time responses use $\rho=5$, $\omega_c=0.2$ and $N\in\{0,3,10\}$, with zero
initial conditions for the plant and all controller states.

\begin{figure*}[t]
\centering
\includegraphics[width=0.9\textwidth]{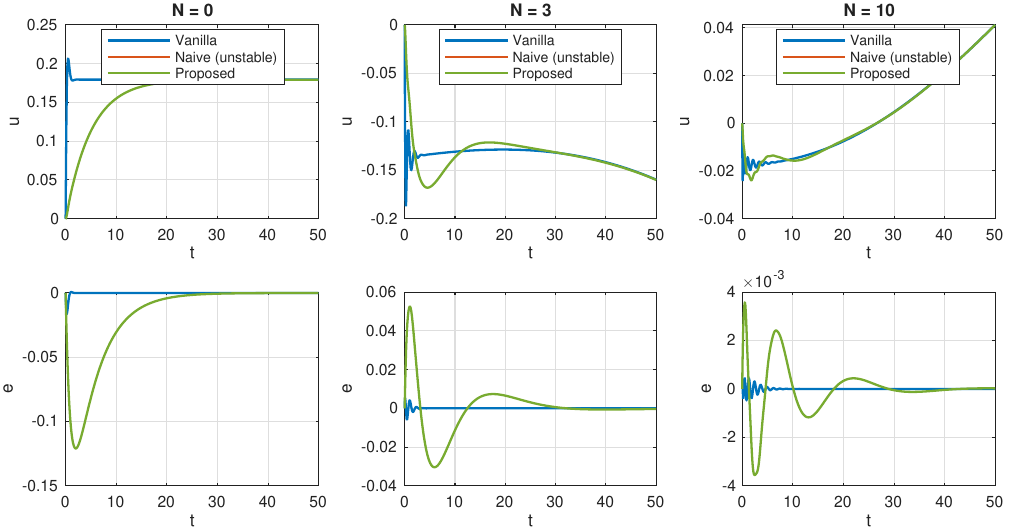}
\caption{Nominal plant ($\Delta A=\Delta B=0$, $\rho=5$, $\omega_c=0.2$).
Control input $u$ (top) and regulated error $e$ (bottom) for the vanilla, naive
and proposed architectures at $N=0,3,10$. The naive series filter is unstable
and is not plotted. The vanilla input peaks and oscillates during the
transient, whereas the proposed input is smooth and regulation $e\to0$ is
preserved. Time window: $t\in[0,50]$\,s.}
\label{fig:nominal}
\end{figure*}

Figure~\ref{fig:nominal} shows the control input (top) and the regulated error
(bottom), one column per order, for the nominal plant. The vanilla controller regulates the error to
zero, but its input shows an initial peak followed by a burst of fast
oscillations during the first seconds of the transient. The proposed
architecture passes this command through $W(s)$ and delivers a markedly
smoother input to the plant. Once the transient is over, the two inputs
coincide, since both converge to the same polynomial steady-state input
$u^\star=-Pw/B$. The price of the smoothing is a slower error transient.
Asymptotic regulation $e\to0$ and closed-loop stability are nonetheless
preserved for every $N$, exactly as Lemma~\ref{lem:convergence} predicts. The naive series filter of Equation~\eqref{eq:naive} is unstable for all three orders (its trajectories
diverge and are omitted from the plots). 

\begin{figure*}[t]
\centering
\includegraphics[width=0.9\textwidth]{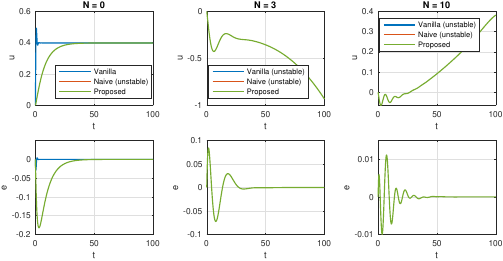}
\caption{Perturbed plant ($\Delta A=0.3$, $\Delta B=-0.55$, $\rho=5$,
$\omega_c=0.2$), with controllers designed on the nominal model. Top: control
input $u$. Bottom: regulated error $e$.
Vanilla is unstable for $N=3,10$ and naive for all $N$ (not plotted). The
proposed scheme satisfies \eqref{eq:sg} and regulates
$e\to0$, as guaranteed by Theorem~\ref{thm:robust}.}
\label{fig:uncertain}
\end{figure*}

Figure~\ref{fig:uncertain} repeats the simulations of
Fig.~\ref{fig:nominal} on the perturbed plant. The proposed architecture
remains stable and regulates the error to zero for every $N$
($|e(100)|<10^{-6}$), with an input as smooth as in the nominal case. For
$N=0$ the steady-state input settles at $u^\star=-Pw/B=0.40$, compared with
$0.18$ in the nominal case. The vanilla loop is unstable for $N=3$ and
$N=10$, and the naive loop for all $N$. 

\begin{figure*}[t]
\centering
\includegraphics[width=0.9\textwidth]{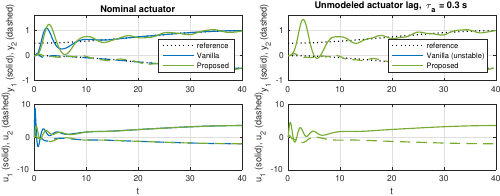}
\caption{Three-mass chain, LQR design with $\mu=10^{-3}$ ($N=3$,
$\omega_c=0.5$). Top: outputs $y_1$ (solid) and $y_2$ (dashed) with their
cubic references (dotted). Bottom: inputs $u_1$ (solid) and $u_2$ (dashed).
Left: nominal actuator. Right: unmodeled first-order actuator lag with time
constant $\tau_a=0.3$\,s, for which the vanilla loop is unstable (not
plotted). The naive loop is unstable in both cases.}
\label{fig:mimo}
\end{figure*}

Next, we examine the benefit of the proposed architecture in improving the robustness of the closed-loop system to input time delays. We define the
\emph{input delay margin} as the largest constant delay $\tau>0$ for which the
closed loop remains stable. We compute it from the loop transfer seen at the
plant input as the phase margin divided by the gain-crossover frequency,
minimized over all crossovers, and report it for the three architectures over
three parameter sweeps: the order $N$ of the internal model (Table~\ref{tab:dm_N},
with $\rho=1$, $\omega_c=0.2$), the closed-loop pole radius $\rho$
(Table~\ref{tab:dm_rho}, with $N=0$, $\omega_c=0.5$), and the filter cutoff
$\omega_c$ (Table~\ref{tab:dm_wc}, with $N=0$, $\rho=1$). 

{In every configuration of
Tables~\ref{tab:dm_N}--\ref{tab:dm_wc} in which the vanilla loop is stable, the
proposed architecture has a larger input delay margin, usually several times
larger. In these examples, filtering the input therefore does more than smooth
the command (Fig.~\ref{fig:nominal}): it also improves delay robustness
relative to the unfiltered design. This improvement is a numerical
observation, and an analytical
characterization of this comparison is left to future work.} The naive series filter is never
competitive. Where it remains stable (small $N$, small $\rho$), its delay
margin is already below the vanilla value, and as the regulator is made more
aggressive, through larger $N$ (Table~\ref{tab:dm_N}) or larger $\rho$
(Table~\ref{tab:dm_rho}), it loses stability.

{A heuristic explanation of the observed improvement is as follows.} The filter $W(s)$ is
low-pass, with $W(0)=I_m$ and $W(j\omega)\to 0$ as $\omega\to\infty$. It enters
in series with the loop seen at the plant input but leaves the spectrum of
$F$, and hence nominal stability, untouched. Attenuating the high-frequency
loop gain lowers the gain-crossover frequency, resulting into a larger tolerable delay. The vanilla law has, on the contrary, crosses over at a higher frequency and is correspondingly more
fragile to delay.

The individual sweeps are consistent with this reading. Increasing $N$
(Table~\ref{tab:dm_N}) makes the regulator more aggressive, which raises the
crossover and shrinks the delay margin of both the vanilla and the proposed
loops, although the proposed loop stays above vanilla throughout. Increasing
$\rho$ (Table~\ref{tab:dm_rho}) shrinks the vanilla margin monotonically, whereas the proposed margin quickly levels off; its crossover is then set by the filter rather than by the
regulator. Decreasing $\omega_c$ (Table~\ref{tab:dm_wc}) makes the filter act
earlier in frequency, which pushes the crossover further down and enlarges the
proposed margin. As $\omega_c$ increases, the filter approaches the identity
and the proposed margin decreases toward the vanilla value.

\begin{table}[htbp]
\centering
\caption{Input delay margin (s) versus the order $N$ of the internal model ($\rho=1$, $\omega_c=0.2$). Perturbed: $\Delta A=0.3$, $\Delta B=-0.55$; unst.: closed loop unstable.}
\label{tab:dm_N}
\footnotesize
\setlength{\tabcolsep}{3.5pt}
\begin{tabular}{ccccccc}
\toprule
 & \multicolumn{3}{c}{Nominal} & \multicolumn{3}{c}{Perturbed} \\
\cmidrule(lr){2-4}\cmidrule(lr){5-7}
$N$ & Van. & Naive & Prop. & Van. & Naive & Prop. \\
\midrule
0 & 1.4827 & 0.5646 & 8.4835 & 2.1691 & 0.8555 & 13.3882 \\
1 & 0.7111 & unst. & 2.6574 & 0.6206 & unst. & 3.7772 \\
2 & 0.5348 & unst. & 1.1983 & 0.1914 & unst. & 2.1700 \\
3 & 0.4328 & unst. & 0.6558 & unst. & unst. & 1.6528 \\
4 & 0.3561 & unst. & 0.4897 & unst. & unst. & 0.4557 \\
5 & 0.2996 & unst. & 0.4120 & unst. & unst. & 0.1608 \\
6 & 0.2574 & unst. & 0.3597 & unst. & unst. & 0.0025 \\
7 & 0.2251 & unst. & 0.3193 & unst. & unst. & unst. \\
8 & 0.1998 & unst. & 0.2867 & unst. & unst. & unst. \\
9 & 0.1794 & unst. & 0.2600 & unst. & unst. & unst. \\
10 & 0.1628 & unst. & 0.2377 & unst. & unst. & unst. \\
\bottomrule
\end{tabular}
\end{table}

\begin{table}[htbp]
\centering
\caption{Input delay margin (s) versus pole radius $\rho$ ($N=0$, $\omega_c=0.5$). Perturbed: $\Delta A=0.3$, $\Delta B=-0.55$; unst.: closed loop unstable.}
\label{tab:dm_rho}
\footnotesize
\setlength{\tabcolsep}{3.5pt}
\begin{tabular}{ccccccc}
\toprule
 & \multicolumn{3}{c}{Nominal} & \multicolumn{3}{c}{Perturbed} \\
\cmidrule(lr){2-4}\cmidrule(lr){5-7}
$\rho$ & Van. & Naive & Prop. & Van. & Naive & Prop. \\
\midrule
1 & 1.4827 & 0.7253 & 3.8756 & 2.1691 & 1.1639 & 6.2397 \\
2 & 0.5285 & 0.0287 & 3.2069 & 0.7496 & 0.0099 & 5.0036 \\
3 & 0.3108 & unst. & 3.1384 & 0.4407 & unst. & 4.8256 \\
4 & 0.2192 & unst. & 3.1294 & 0.3107 & unst. & 4.7747 \\
5 & 0.1691 & unst. & 3.1297 & 0.2397 & unst. & 4.7544 \\
6 & 0.1376 & unst. & 3.1314 & 0.1950 & unst. & 4.7444 \\
7 & 0.1160 & unst. & 3.1331 & 0.1643 & unst. & 4.7390 \\
8 & 0.1002 & unst. & 3.1345 & 0.1419 & unst. & 4.7356 \\
9 & 0.0882 & unst. & 3.1356 & 0.1249 & unst. & 4.7335 \\
10 & 0.0788 & unst. & 3.1365 & 0.1115 & unst. & 4.7320 \\
\bottomrule
\end{tabular}
\end{table}

\begin{table}[htbp]
\centering
\caption{Input delay margin (s) versus filter cutoff $\omega_c$ ($N=0$, $\rho=1$). Perturbed: $\Delta A=0.3$, $\Delta B=-0.55$; unst.: closed loop unstable.}
\label{tab:dm_wc}
\footnotesize
\setlength{\tabcolsep}{3.5pt}
\begin{tabular}{ccccccc}
\toprule
 & \multicolumn{3}{c}{Nominal} & \multicolumn{3}{c}{Perturbed} \\
\cmidrule(lr){2-4}\cmidrule(lr){5-7}
$\omega_c$ & Van. & Naive & Prop. & Van. & Naive & Prop. \\
\midrule
0.1 & 1.4827 & 0.4955 & 16.3068 & 2.1691 & 0.7144 & 25.5311 \\
0.2 & 1.4827 & 0.5646 & 8.4835 & 2.1691 & 0.8555 & 13.3882 \\
0.3 & 1.4827 & 0.6247 & 5.8991 & 2.1691 & 0.9742 & 9.3830 \\
0.4 & 1.4827 & 0.6778 & 4.6254 & 2.1691 & 1.0759 & 7.4072 \\
0.5 & 1.4827 & 0.7253 & 3.8756 & 2.1691 & 1.1639 & 6.2397 \\
0.6 & 1.4827 & 0.7681 & 3.3867 & 2.1691 & 1.2408 & 5.4738 \\
0.7 & 1.4827 & 0.8070 & 3.0459 & 2.1691 & 1.3082 & 4.9356 \\
0.8 & 1.4827 & 0.8423 & 2.7968 & 2.1691 & 1.3677 & 4.5386 \\
0.9 & 1.4827 & 0.8747 & 2.6080 & 2.1691 & 1.4204 & 4.2347 \\
1.0 & 1.4827 & 0.9043 & 2.4609 & 2.1691 & 1.4673 & 3.9954 \\
\bottomrule
\end{tabular}
\end{table}

\subsection{A multivariable example}\label{sec:num_mimo}

Consider a chain of three unit masses connected to each other and to two
walls by four unit springs, with dampers of coefficient $0.02$ ($n=6$). With
state $x=[\,q^\top\ \dot q^\top\,]^\top$, where $q\in\mathbb{R}^3$ collects the
mass positions, the plant matrices in \eqref{eq:plant}--\eqref{eq:err} are
\begin{equation}\label{eq:mimo_plant}
A=\begin{bmatrix} 0_3 & I_3\\ -K_s & -0.02\,K_s\end{bmatrix},\qquad
B=\begin{bmatrix} 0_{3\times2}\\ [\,\varepsilon_1\ \ \varepsilon_3\,]\end{bmatrix},
\end{equation}
and $C=[\,\varepsilon_2\ \ \varepsilon_3\,]^\top[\,I_3\ \ 0_3\,]\in\mathbb{R}^{2\times6}$,
where $\varepsilon_i$ is the $i$th column of $I_3$ and $K_s\in\mathbb{R}^{3\times3}$ is the tridiagonal stiffness matrix
with $2$ on the diagonal and $-1$ on the first sub- and superdiagonals. Forces act on the first and third mass, and the outputs are the
positions of the second and third mass, $y=Cx\in\mathbb{R}^2$, so $m=p=2$ and
one output is not collocated with the inputs. 

Each output tracks a cubic reference $r(t)\in\mathbb{R}^2$ ($N=3$, so
$w(t)\in\mathbb{R}^4$ and $n_\eta=8$) over $T=40$\,s. With
$\theta=t/T$,
\[
r_1(t)=0.5+1.5\,\theta^2-\theta^3,\qquad r_2(t)=-1.5\,\theta^2+\theta^3,
\]
so that $y_1$ moves from $0.5$ to $1$ and $y_2$ from $0$ to $-0.5$. A constant force
of $0.2$ acts on the unactuated mass. The internal model is that of Section~\ref{sec:num_setup}, and
the filter applies \eqref{eq:num_filter} to each input channel,
$W(s)=\big(1-(s/(s+\omega_c))^{N+1}\big)I_2$, with $\omega_c=0.5$. For the gains $K_x,K_\eta$,
LQR is the natural choice in the multivariable case. We compute them for
$(A,B)$ augmented with the internal model, with state weight
$\mathrm{blkdiag}(10\,C^\top C+0.1\,I_6,\,I_8)$ and input weight $\mu\,I_2$,
so that decreasing $\mu$ makes the regulator more aggressive.

Table~\ref{tab:mimo} compares three such designs. As $\mu$ decreases, the
vanilla controller demands larger inputs and becomes more fragile. Its peak input more than doubles, the largest unmodeled first-order
actuator lag it tolerates drops from $0.175$\,s to $0.100$\,s, and its delay
margin for a common delay in both input channels shrinks by a factor of three.
The proposed controller is essentially insensitive to $\mu$. Its peak input
stays near $3.6$, it tolerates actuator lags longer than $1$\,s, and its delay
margin stays between $0.44$ and $0.59$\,s, which is $4$ to $19$ times the
vanilla value.
Figure~\ref{fig:mimo} shows the time responses of the design with
$\mu=10^{-3}$. With the nominal actuator, both controllers track the
references, the proposed one with a slower transient and a peak input of
$3.59$ instead of $8.95$. An unmodeled actuator lag of $0.3$\,s, a
perturbation outside the setting of Theorem~\ref{thm:robust}, separates the two
designs: the vanilla loop becomes unstable, whereas the proposed loop still
tracks both references.

\begin{table}[htbp]
\centering
\caption{Three-mass chain ($N=3$, $\omega_c=0.5$) for increasingly aggressive LQR designs (input weight $\mu I_2$): peak input on the tracking task, largest tolerable time constant $\tau_{\max}$ of an unmodeled actuator lag, and input delay margin (s). The naive series filter is unstable for all three designs.}
\label{tab:mimo}
\footnotesize
\setlength{\tabcolsep}{3.5pt}
\begin{tabular}{ccccccc}
\toprule
 & \multicolumn{2}{c}{$\max|u|$} & \multicolumn{2}{c}{$\tau_{\max}$ (s)} & \multicolumn{2}{c}{Delay margin (s)} \\
\cmidrule(lr){2-3}\cmidrule(lr){4-5}\cmidrule(lr){6-7}
$\mu$ & Van. & Prop. & Van. & Prop. & Van. & Prop. \\
\midrule
$10^{-2}$ & 5.40 & 3.59 & 0.175 & 1.059 & 0.102 & 0.439 \\
$10^{-3}$ & 8.95 & 3.59 & 0.127 & 1.117 & 0.060 & 0.581 \\
$10^{-4}$ & 12.61 & 3.65 & 0.100 & 1.132 & 0.031 & 0.590 \\
\bottomrule
\end{tabular}
\end{table}

\section{Conclusions}\label{sec:conclusions}

This paper addressed the transient input peaking of internal model
principle-based controllers for polynomial exogenous signals. In the proposed
architecture, the internal model regulates the error of an ideal, unfiltered
copy of the plant, the physical plant receives a low-pass filtered command,
and an auxiliary plant copy keeps the filter out of the loop that the
regulator stabilizes. Under blocking conditions on the filter at the origin,
the regulator and the filter can therefore be designed independently without
losing internal stability or regulation. When the plant is uncertain, a
small-gain condition on explicit loop maps preserves both properties. In the numerical
examples, the proposed controller delivered a smoother input and larger
input-delay margins than the unfiltered one, and its tolerance to unmodeled
actuator dynamics did not degrade as the regulator was made more aggressive.
The naive series filter was unstable in almost every configuration.

Open questions include an analytical characterization of the
delay-margin improvement,
explicit bounds on the error transient,
optimization-based filter design, and extensions to nonlinear plants, input
saturation, and piecewise-polynomial exogenous signals.

\bibliographystyle{IEEEtran}
\bibliography{refs}

\end{document}